\documentclass[11pt]{article}
\usepackage{amssymb,amsmath}
\newcommand{\be}{\begin{equation}}
\newcommand{\ee}{\end{equation}}
\newcommand{\bea}{\begin{eqnarray}}
\newcommand{\eea}{\end{eqnarray}}
\newcommand{\mbf}[1]{\mbox{\boldmath$#1$}}

\begin{document}

\title{Complete Bell polynomials and new generalized identities for polynomials of higher order}
\author{Boris Rubinstein,
\\Stowers Institute for Medical Research
\\1000 50$^{}\mbox{th}$ St., Kansas City, MO 64110, U.S.A.}
\date{\today}

\maketitle
\begin{abstract}
The relations between the Bernoulli and Eulerian polynomials of higher order and the complete Bell polynomials are found that lead to new identities for the Bernoulli and Eulerian polynomials and numbers of higher order. General form of these identities is considered and generating function for polynomials satisfying this general identity is found.
\end{abstract}

\section{Introduction}
The history of the Bernoulli polynomials $B_n(s)$
counts more than 250 years, as
L. Euler first studied them for arbitrary values of the argument.
He introduced the Euler polynomials $E_n(s)$,
the generalization of these
polynomials to the so-called Eulerian polynomials
$H_n(s,\rho)$ was made by Frobenius \cite{Frobenius}.
The Bell polynomials \cite{Comtet1974} also appeared to be useful in combinatorial applications.


In 1920s N. N\"orlund \cite{NorlundMemo}
introduced the
Bernoulli $B^{(m)}_n(s|{\bf d})$
and Euler polynomials $E^{(m)}_n(s|{\bf d})$
of higher order adding parameters ${\bf d}=\{d_i\}, \ 1\le i \le m$.
Similar extension for the Eulerian polynomials
$H^{(m)}_n(s,\mbf{\rho}|{\bf d})$ was made by
L. Carlitz in \cite{Carlitz1960}.

The Bernoulli and Eulerian polynomials of higher order
appear to be useful for the description of the
restricted partition function 
defined as a number of integer nonnegative solutions to
Diophantine equation ${\bf d}\cdot{\bf x}=s$
(see \cite{Rama03}).
The author showed in \cite{Rub04}
that the restricted partition function can be expressed through the
Bernoulli polynomials of higher order only. 

It appears that there exist relations that connect the Bernoulli, Euler and related polynomials and their higher order generalizations 
to the complete Bell polynomials. These relations presented in Section 3 give rise to new identities for the polynomials of higher order, and the identity for the Eulerian polynomials is discussed in Section 4. In Section 5 we introduce a more general form of the abovementioned identities and find the generating functions for the polynomials satisfying these relations.
The Section 6 is devoted to derivation of new relations for the 
Bernoulli polynomials of higher order.

\section{Bernoulli and Eulerian polynomials and their generalization to higher order}
The Bernoulli $B_n(s)$ and Euler $E_n(s)$
polynomials are defined through the corresponding generating functions 
\be
\frac{te^{st}}{e^{t}-1} =
\sum_{n=0}^{\infty} B_n(s)
\frac{t^{n}}{n!}, \ \ \
\frac{2e^{st}}{e^{t}+1} =
\sum_{n=0}^{\infty} E_n(s)
\frac{t^{n}}{n!}. 
\label{RegPolyGF}
\ee
Frobenius \cite{Frobenius} studied the
so-called Eulerian polynomials
$H_n(s,\rho)$ satisfying the generating function
\be
e^{st}\frac{1-\rho}{e^t-\rho} = \sum_{n=0}^{\infty} H_n(s,\rho)
\frac{t^n}{n!}, \ \ (\rho \ne 1),
\label{EulerRegGFnew}
\ee
which reduces to definition of the Euler polynomials at fixed value of the
parameter $E_n(s) = H_n(s,-1).$
The Apostol-Bernoulli \cite{Apostol1951} and Apostol-Euler \cite{Luo2005} are the particular cases of the Eulerian polynomials.

Introduce the generalized Eulerian polynomials through the generating function
\be
e^{st} \frac{1-\alpha\rho}{e^t-\rho} =
\sum_{n=0}^{\infty} H_n(s,\rho,\alpha) \frac{t^{n}}{n!}\;.
\label{EulerianGenGF_basic}
\ee
It is easy to see that
$$
B_n(s) = n H_{n-1}(s,1,0), \
E_n(s) = H_n(s,-1,1), \
H_n(s,\rho) = H_n(s,\rho,1).
$$

Introduce the generalized Eulerian polynomials of higher order
$H^{(m)}_n(s,\mbf{\rho},\mbf{\alpha}|{\bf d})$ as follows
\be
e^{st} \prod_{i=1}^m \frac{1-\alpha_i\rho_i}{e^{d_i t}-\rho_i} =
\sum_{n=0}^{\infty} H^{(m)}_n(s,\mbf{\rho},\mbf{\alpha}|{\bf d})
\frac{t^{n}}{n!}\;.
\label{EulerianGenGF}
\ee
The generating function for the Bernoulli polynomials of higher
order $B^{(m)}_n(s|{\bf d})$ reads \cite{bat53}:
\be
e^{st} \prod_{i=1}^m \frac{d_i t}{e^{d_it}-1} =
\sum_{n=0}^{\infty} B^{(m)}_n(s|{\bf d})
\frac{t^{n}}{n!}\;.
\label{BernoulliGF}
\ee
It can be checked by comparison that 
$$
B^{(m)}_n(s|{\bf d})=
\frac{n!\pi_1}{(n-m)!} H^{(m)}_{n-m}(s,{\bf 1},{\bf 0}|{\bf d}), \ \
\pi_1 = \prod_{1=1}^m d_i\;.
$$
It is easy to see that the Euler $E^{(m)}_n(s|{\bf d})$ \cite{bat53} and Eulerian polynomials of higher order $H^{(m)}_n(s,\mbf{\rho}|{\bf d})$ \cite{Carlitz1960}
can be written as 
\bea
E^{(m)}_n(s|{\bf d})&=&H^{(m)}_n(s,-{\bf 1},{\bf 1}|{\bf d}),
\label{EulerianGenCases} \\
H^{(m)}_n(s,\mbf{\rho}|{\bf d})&=&H^{(m)}_n(s,\mbf{\rho},{\bf 1}|{\bf d}). \nonumber
\eea

\section{Bernoulli and Eulerian polynomials of higher order
as Bell polynomials}

The complete Bell polynomials ${\bf B}_n(a_1,a_2,\ldots)$ \cite{Riordan} are defined through the generating function
\be
\exp \left(\sum_{i=1}^{\infty}\frac{a_i}{i!}t^i \right) =
\sum_{i=0}^{\infty} \frac{{\bf B}_i(a_1,a_2,\ldots)}{i!}t^i.
\label{BellGF}
\ee
Consider Taylor series
$$
\ln \frac{dt}{e^{dt}-1} =
\sum_{i=1}^{\infty} (-1)^{i-1} \frac{d^i B_i}{i} \frac{t^i}{i!},
$$
where $B_i$ is the Bernoulli number. Using the above expansion in
(\ref{BernoulliGF}) we have
\be
e^{st}
\exp\left(\sum_{i=1}^{\infty} (-1)^{i-1}\frac{B_i \sigma_i}{i}
\frac{t^{i}}{i!} \right)\; =
\sum_{i=0}^{\infty} {\bf B}_i(s+a_1,a_2,\ldots)\frac{t^{i}}{i!} ,
\label{BernoulliGF1}
\ee
where
$a_i = (-1)^{i-1}B_i \sigma_i/i$ and $\sigma_k=\sigma_k({\bf d})=\sum_{i=1}^{m}d_i^k$ denotes a power sum, so that
\be
B^{(m)}_n(s|{\bf d}) = {\bf B}_n(s+a_1,a_2,\ldots).
\label{BernoulliBell}
\ee
Setting here $m=1,\ d=1$ we find 
$$
B_n(s) = {\bf B}_n(s+a_1,a_2,\ldots), \ \ \
a_i = (-1)^{i-1}\frac{B_i}{i}.
$$
Consider Taylor series for the Eulerian polynomials
$$
\ln \frac{1-\rho}{e^{d t}-\rho} = -\frac{dt}{1-\rho}
-\sum_{i=2}^{\infty} \frac{\rho d^i H_{i-1}(\rho)}{1-\rho}\frac{t^i}{i!}.
$$
where $H_{i}(\rho)=H_{i}(0,\rho)$ is the Eulerian number.
In a particular case with $\rho_i = \rho$  in
(\ref{EulerianGenGF})
using the above expansion we have
\be
\exp\left(st-\frac{\sigma_1t}{1-\rho}\right)
\exp\left(-\sum_{i=2}^{\infty} \frac{\rho H_{i-1}(\rho) \sigma_i}{1-\rho}
\frac{t^{i}}{i!} \right)\; =
\sum_{i=0}^{\infty}
{\bf B}_i(s-\sigma_1+a_1,a_2,\ldots)\frac{t^{i}}{i!},
\label{EulerianGF1}
\ee
where
$a_i = -\rho H_{i-1}(\rho) \sigma_i/(1-\rho)$.
Comparing (\ref{EulerianGF1}) to (\ref{BellGF}) we obtain
\be
H^{(m)}_n(s,\rho|{\bf d}) = {\bf B}_n(s-\sigma_1+a_1,a_2,\ldots).
\label{EulerianBell}
\ee
The relations 
(\ref{BernoulliBell},\ref{EulerianBell}) 
are the particular cases of a general identity
\be
H^{(m)}_n(s,\rho,\alpha|{\bf d}) =
{\bf B}_n(s-\sigma_1+a_1,a_2,\ldots)
\label{EulerianGenBell}
\ee
with $a_i = -\rho H_{i-1}(\rho,\alpha) \sigma_i/(1-\alpha\rho)$ and
$H_{i}(\rho,\alpha) = H_{i}(0,\rho,\alpha)$.
Note that the value of the polynomial is equal to the complete Bell polynomial of the same order with parameters proportional to the value of the corresponding polynomial at zero.

\section{Generalized Eulerian polynomials identities}

Consider a normalized distribution $P(t)$ that can be characterized
by the sequence of moments $\mu_n$ or cumulants $\kappa_n$. The
moments generation function is $P(t)$ itself and cumulants
generation function is $\ln P(t)$, i.e.,
\be P(t) =
1+\sum_{n=1}^{\infty} \frac{\mu_n}{n!} t^n, \ \ \ \ln P(t) =
\sum_{n=1}^{\infty} \frac{\kappa_n}{n!} t^n.
\label{moment_cumulants_GF}
\ee
It is worth to note that in the theory of symmetric functions $\mu_n$ play role of the complete homogeneous symmetric polynomials, while $\kappa_n$ stand for the symmetric power sums.
The above definitions imply the relation between the moments and
cumulants
\be \mu_n = \kappa_n+\sum_{k=1}^{n-1} \binom{n-1}{k-1}
\kappa_k \mu_{n-k}.
\label{moments_cumulant_rel}
\ee
On the other
hand it follows from (\ref{moment_cumulants_GF}) and (\ref{BellGF})
that
\be \mu_n = {\bf B}_n(\kappa_1,\kappa_2,\ldots),
\label{moment_cumulants_Bell}
\ee
so that
\be
{\bf B}_n(\kappa_1,\kappa_2,\ldots) =\kappa_n+ \sum_{k=1}^{n-1}
\binom{n-1}{k-1} \kappa_k {\bf B}_{n-k}(\kappa_1,\kappa_2,\ldots).
\label{moments_cumulant_relBell}
\ee
Using the relation (\ref{EulerianGenBell}) and the definition (\ref{BellGF}) one finds:
\be
\ln \sum_{n=0}^{\infty} \frac{H^{(m)}_n(s,\rho,\alpha|{\bf d})}{n!}t^n =
(s-\sigma_1)t-\frac{\rho}{1-\alpha\rho}
\sum_{i=1}^{\infty} \frac{H_{i-1}(0,\rho,\alpha)\sigma_i}{i!}\;t^i.
\label{EulerianGenHO_mom_cum}
\ee
The last relation leads to the following identity:
$$
H^{(m)}_n(s,\rho,\alpha|{\bf d}) =
(s-\sigma_1)H^{(m)}_{n-1}(s,\rho,\alpha|{\bf d})
-\frac{\rho}{1-\alpha\rho}\sum_{k=1}^{n}
\binom{n-1}{k-1} H_{k-1}(0,\rho,\alpha)\sigma_k
H^{(m)}_{n-k}(s,\rho,\alpha|{\bf d}),
$$
that is rewritten as
\be
H^{(m)}_{n+1}(s,\rho,\alpha|{\bf d}) =
(s-\sigma_1)H^{(m)}_{n}(s,\rho,\alpha|{\bf d})-
\frac{\rho}{1-\alpha\rho}\sum_{k=0}^{n}
\binom{n}{k} H_{k}(0,\rho,\alpha)\sigma_{k+1}
H^{(m)}_{n-k}(s,\rho,\alpha|{\bf d}).
\label{EulerianGenHO_Euler_relation}
\ee
Setting here $m=1, \; d_1=1$ we find for the generalized Eulerian polynomials
\be
H_{n+1}(s,\rho,\alpha) =
(s-1)H_{n}(s,\rho,\alpha)-
\frac{\rho}{1-\alpha\rho}\sum_{k=0}^{n}
\binom{n}{k} H_{k}(0,\rho,\alpha) H_{n-k}(s,\rho,\alpha).
\label{EulerianGenPoly_relation}
\ee

\section{General case}
Consider a more general problem of construction of polynomials
$P_n(s,r)$ that can be represented by the Bell complete polynomials with
parameters $a_i$ being proportional to the values of the same polynomial at some $s=r$, that satisfy the generalized version of (\ref{EulerianGenHO_mom_cum})
\be
\ln \sum_{n=0}^{\infty} \frac{P_n(s,r)}{n!}t^n =
Q_1(s,r)t+Q_2(s,r)
\sum_{i=1}^{\infty} \frac{P_{i-1}(r,r)}{i!}\;t^i,
\label{polynomial_HO_mom_cum}
\ee
where $Q_1(s,r)$ and $Q_2(s,r)$ are polynomials in $s$.
Then we find
\be
P_{n+1}(s,r) =
Q_1(s,r)P_{n}(s,r)+
Q_2(s,r)\sum_{k=0}^{n}
\binom{n}{k} P_{k}(r,r) P_{n-k}(s,r).
\label{polynomial_Poly_relation}
\ee
We want to find the condition on the generating function $G(s,r,t)$ for the polynomials $P_n(s,r)$. Multiplying the equation (\ref{polynomial_Poly_relation}) by $t^n/n!$ and summing over $n$ we obtain
$$
\sum_{n=0}^{\infty} \frac{P_{n+1}(s,r)}{n!} t^n =
\sum_{n=0}^{\infty} \left[
Q_1(s,r)P_{n}(s,r)+
Q_2(s,r)\sum_{k=0}^{n}
\binom{n}{k} P_{k}(r,r) P_{n-k}(s,r)
\right]\frac{t^n}{n!}.
$$
The r.h.s. of the above equation evaluates to
$Q_1(s,r)G(s,r,t)+Q_2(s,r) G(r,r,t)G(s,r,t)$, while the l.h.s. is obtained by differentiation of $G(s,r,t)$ w.r.t. $t$, so that
\be
\frac{\partial G(s,r,t)}{\partial t} = 
Q_1(s,r)G(s,r,t)+Q_2(s,r) G(r,r,t)G(s,r,t).
\label{GF_diff_eq}
\ee

Consider a particular case of (\ref{GF_diff_eq}) in the form 
\be
\frac{\partial G(s,r,t)}{\partial t} = Q_1(s)G(s,r,t)+Q_2(r)
G(r,r,t)G(s,r,t). 
\label{GF_diff_eq_gen} 
\ee 
The polynomials $P_{n}(s,r)$ generated by $G(s,r,t)$ satisfy the identity
\be
P_{n+1}(s,r)=Q_1(s)P_{n}(s,r)+Q_2(r)\sum_{k=0}^{n}\binom{n}{k}P_{k}(r,r)
P_{n-k}(s,r).
\label{z2} 
\ee
We look for solution in
the form $G(s,r,t)=H(s,t)F(r,t)$ that leads to the equation 
\be
\frac{\partial F(r,t)}{\partial t} =\left[Q_1(s)-\frac{\partial \ln
H(s,t)}{\partial t}\right]F(r,t) + Q_2(r)H(r,t) F^2(r,t).
\label{diff_eqF} 
\ee 
The term in the square brackets in the last equation 
should be independent of $s$, so that the following condition holds
$$
Q_1(s)-\frac{\partial \ln H(s,t)}{\partial t} = M'(t)
$$
for some function $M(t)$ defined up to arbitrary constant. The
condition leads to relation 
\be H(s,t) = e^{Q_1(s)t-M(t)}\;,
\label{reh_H} 
\ee 
and the equation (\ref{diff_eqF}) reduces to the Bernoulli equation 
\be 
\frac{\partial F(r,t)}{\partial t}
=M'(t)F(r,t) + Q_2(r)e^{Q_1(r)t-M(t)} F^2(r,t),
\label{diff_eqF_Bernoulli} 
\ee 
which has the solution
\be
F(r,t) = \frac{Q_1(r)e^{M(t)}}{C Q_1(r) - Q_2(r) e^{Q_1(r)t}}.
\label{a1}
\ee
The solution of (\ref{GF_diff_eq_gen}) reads
\be
G(s,r,t)=H(s,t)F(r,t) = 
\frac{Q_1(r)e^{Q_1(s)t}}{C Q_1(r) - Q_2(r)e^{Q_1(r)t}},
\label{a2}
\ee
where the integration constant $C$ is determined from
the condition
\be
G(s,r,0) = G_0(s,r) = \frac{Q_1(r)}{C Q_1(r) - Q_2(r)}.
\label{a3}
\ee
Finally we obtain 
\be G(s,r,t) =\frac{G_0(s,r)Q_1(r) e^{Q_1(s)t}}
{Q_1(r)+G_0(s,r)Q_2(r)(1-e^{Q_1(r)t})}\;. 
\label{GF_diff_eq_gen_sol}
\ee
Setting in the last equation $r=0,\ Q_1(s)=s-1, \ Q_2(0)=-\rho/(1-\alpha\rho)$ and $G_0(s,0) = (1-\alpha\rho)/(1-\rho)$ we obtain the generating function in the l.h.s. of (\ref{EulerianGenGF_basic}).

It is not difficult to generalize the above generating function to the case of polynomials of higher order. We define the generating function $G(s,r,t,{\bf d})$ for the polynomials of higher order $P^{(m)}_n(s,r|{\bf d})$ as follows
\be
G(s,r,t,{\bf d})=e^{[Q_1(s)-Q_1(0)]t}\prod_{i=1}^{m}
\frac{e^{Q_1(0) d_1 t}}
{G_0^{-1}(s,r)+Q_2(r)(1-e^{Q_1(r)d_i t})/Q_1(r)}\;.
\label{polynomial_PolyGenGF_basic}
\ee
The polynomial identity reads:
\be
P^{(m)}_{n+1}(s,r|{\bf d}) =
(Q_1(s)-Q_1(0)(1-\sigma_1))P^{(m)}_{n}(s,r|{\bf d})+
Q_2(r)\sum_{k=0}^{n}
\binom{n}{k} P_{k}(r,r)\sigma_{k+1} P^{(m)}_{n-k}(s,r|{\bf d}).
\label{polynomial_PolyHO_Euler_relation}
\ee
The case $r=s$ leads to the general solution 
\be G(s,t)
=\frac{G_0(s)Q_1(s) e^{Q_1(s)t}}
{Q_1(s)+G_0(s)Q_2(s)(1-e^{Q_1(s)t})}\;.
\label{GF_diff_eq_gen_sol_r=s} 
\ee
It should be noted that the solutions
(\ref{GF_diff_eq_gen_sol},\ref{GF_diff_eq_gen_sol_r=s}) are valid
only when $Q_1(s)$ is identically nonzero. Otherwise the equation (\ref{GF_diff_eq})
reduces to 
\be \frac{\partial G(s,r,t)}{\partial t} = Q_2(s,r)
G(r,r,t)G(s,r,t), 
\label{GF_diff_eq_gen2} 
\ee 
that can be written in the form
$$
\frac{\partial F(r,t)}{\partial t}=G(r,r,t), \ \ \ F(r,t) =
\frac{\ln G(s,r,t)}{Q_2(s,r)}.
$$
From the definition of $F(r,t)$ we find
$G(s,r,t)=\exp[Q_2(s,r)F(r,t)]$ and obtain
$G(r,r,t)=\exp[Q_2(r,r)F(r,t)]$. 
Thus we arrive at the following
equation for $F(r,t)$
$$
\frac{\partial F(r,t)}{\partial t} = \exp(Q_2(r,r)F(r,t)),
$$
which has the solution
$$
F(r,t) = -\frac{\ln Q_2(r,r)(C-t)}{Q_2(r,r)},
$$
and we obtain the solution of equation (\ref{GF_diff_eq_gen2}) 
\be G(s,r,t) =
\exp\left[-\frac{Q_2(s,r)}{Q_2(r,r)}\ln Q_2(r,r)(C-t) \right] =
\left[ Q_2(r,r)(C-t) \right]^{-Q_2(s,r)/Q_2(r,r)}.
\label{GF_diff_eq_gen2_sol} 
\ee 
Setting in
(\ref{GF_diff_eq_gen2_sol}) $t=0$ we obtain the integration constant
$$
C = \frac{1}{Q_2(r,r)}G_0^{-Q_2(r,r)/Q_2(s,r)}\;, \ \ G_0 =
G(s,r,0).
$$
Finally, we arrive at the solution 
\be 
G(s,r,t) = G_0\left[
1-Q_2(r,r)tG_0^{Q_2(r,r)/Q_2(s,r)} \right]^{-Q_2(s,r)/Q_2(r,r)}.
\label{GF_diff_eq_gen2_sol_fin} 
\ee
For $r=s$ equation (\ref{GF_diff_eq_gen2}) reduces to 
\be 
\frac{\partial
G(s,t)}{\partial t} = Q_2(s) G^2(s,t). 
\label{GF_diff_eq_gen2_r=s}
\ee 
The solution in this case is found from 
(\ref{GF_diff_eq_gen_sol_r=s}) in the limit $Q_1(s)\to 0$
\be 
G(s,t) = \frac{G_0(s)}{1-G_0(s)Q_2(s)t}.
\label{GF_diff_eq_gen2_r=s_sol}
\ee

\section{Identities for Bernoulli polynomials and numbers}

Using the relation (\ref{BernoulliBell}) and the definition (\ref{BellGF}) one finds:
\be
\ln \sum_{n=0}^{\infty} \frac{B^{(m)}_n(s|{\bf d})}{n!}t^n =
st+
\sum_{i=1}^{\infty} (-1)^{i-1}\frac{B_i\sigma_i}{i!\;i}\;t^i
.
\label{BernoulliHO_mom_cum}
\ee
The last relation leads to the following identity:
\be
B^{(m)}_n(s|{\bf d}) = sB^{(m)}_{n-1}(s|{\bf d})-
\frac{1}{n}\sum_{k=1}^{n}
\binom{n}{k} (-1)^k B_k\sigma_k B^{(m)}_{n-k}(s|{\bf d}).
\label{BernoulliHO_Euler_relation}
\ee
Setting here $m=1, \; d_1=1$ we find for the Bernoulli polynomials
\be
B_n(s) = sB_{n-1}(s)-
\frac{1}{n}\sum_{k=1}^{n}
\binom{n}{k} (-1)^k B_k B_{n-k}(s).
\label{Bernoulli_Euler_relation}
\ee
Finally setting in the above relation $s=0$ we arrive at the identity for the Bernoulli numbers
\be
\sum_{k=1}^{n}
\binom{n}{k} (-1)^k B_k B_{n-k} = -n B_n.
\label{BernoulliNo_Euler_relation}
\ee
As all odd Bernoulli numbers (except $B_1$) vanish, one can find
\be
\sum_{k=1}^{n}
\binom{n}{k}[1-(-1)^k]B_k\sigma_k B^{(m)}_{n-k}(s|{\bf d}) =
-n \sigma_1 B^{(m)}_{n-1}(s|{\bf d}),
\label{Bernoulli_Euler_relation3}
\ee
from what it follows that
\bea
B^{(m)}_n(s|{\bf d}) &= &(s-\sigma_1)B^{(m)}_{n-1}(s|{\bf d})-
\frac{1}{n}\sum_{k=1}^{n}
\binom{n}{k}B_k\sigma_k B^{(m)}_{n-k}(s|{\bf d}),
\label{Bernoulli_Euler_relation1} \\
B_n(s) & = & (s-1)B_{n-1}(s)-
\frac{1}{n}\sum_{k=1}^{n}
\binom{n}{k} B_k B_{n-k}(s),
\label{Bernoulli_Euler_relation2} \\
B_n & = &-B_{n-1} - \frac{1}{n}\sum_{k=1}^{n}
\binom{n}{k} B_k B_{n-k},
\label{Euler_relation}
\eea
where the relation (\ref{Euler_relation}) is the well-known Euler identity and (\ref{Bernoulli_Euler_relation2}) is a particular case of the sum identity for the Bernoulli polynomials
$$
\sum_{k=0}^n \binom{n}{k} B_k(s) B_{n-k}(t) =
n(s+t-1) B_{n-1}(s+t)-(n-1)B_n(s+t).
$$

Rewrite the identity (\ref{Bernoulli_Euler_relation1}) using the definition of the power sum as follows:
\bea
B^{(m)}_n(s|{\bf d}) &= &(s-\sigma_1)B^{(m)}_{n-1}(s|{\bf d})-
\frac{1}{n}\sum_{i=1}^m\sum_{k=1}^{n}
\binom{n}{k}B_k d_i^k B^{(m)}_{n-k}(s|{\bf d})
\nonumber \\
&= &(s-\sigma_1)B^{(m)}_{n-1}(s|{\bf d})-
\frac{1}{n}\sum_{i=1}^m
\left[-B^{(m)}_{n}(s|{\bf d})+
\sum_{k=0}^{n}
\binom{n}{k}B_k\sigma_k B^{(m)}_{n-k}(s|{\bf d})\right]
\nonumber \\
&= &(s-\sigma_1)B^{(m)}_{n-1}(s|{\bf d})+
\frac{m}{n}B^{(m)}_{n}(s|{\bf d})-
\frac{1}{n}\sum_{i=1}^m
B^{(m+1)}_{n}(s|{\bf d}_i),
\label{Bernoulli_Euler_relation1_upup}
\eea
where
$$
{\bf d}_i = {\bf d} \cup d_i.
$$
The relation (\ref{Bernoulli_Euler_relation1_upup}) produces two new identities
\be
\sum_{i=1}^m
B^{(m+1)}_{n}(s|{\bf d}_i) =
n(s-\sigma_1)B^{(m)}_{n-1}(s|{\bf d})-
(n-m)B^{(m)}_n(s|{\bf d}),
\label{upup_gen}
\ee
and
\be
\sum_{i=1}^m
B^{(m+1)}_{m}(s|{\bf d}_i) =
m(s-\sigma_1) B^{(m)}_{m-1}(s|{\bf d}).
\label{upup_m=n}
\ee
Recalling the identity for the Bernoulli polynomials of higher order
\cite{NorlundMemo}:
$$
B^{(m+1)}_{n}(s+d_i|{\bf d}_i)-B^{(m+1)}_{n}(s|{\bf d}_i) =
B^{(m)}_{n-1}(s|{\bf d}),
$$
we have
\be
\sum_{i=1}^m
B^{(m+1)}_{n}(s+d_i|{\bf d}_i) =
ns B^{(m)}_{n-1}(s|{\bf d}) -(n-m)B^{(m)}_{n}(s|{\bf d}).
\label{upup2}
\ee

Consider a particular case of all $d_i=1$, i.e., ${\bf d} = {\bf 1}$; noting that
$\sigma_1 = m$ we find from (\ref{upup_gen})
the identity derived by N\"orlund in \cite{NorlundMemo}
\be
m B^{(m+1)}_{n}(s|{\bf 1}) =
(m-n)B^{(m)}_{n}(s|{\bf 1}) +
n(s-m)B^{(m)}_{n-1}(s|{\bf 1}).
\label{Bernoulli_Euler_m_plus_1}
\ee
It is easy to show that for arbitrary nonzero $\lambda$ the following relation is valid:
\be
m B^{(m+1)}_{n}(s|\lambda\cdot{\bf 1}) =
(m-n)B^{(m)}_{n}(s|\lambda\cdot{\bf 1}) +
n(s-\lambda m)B^{(m)}_{n-1}(s|\lambda\cdot{\bf 1}).
\label{Bernoulli_Euler_m_plus_1lambda}
\ee
From the last relation for $m=n$ it follows that
\be
B^{(m+1)}_{m}(s|\lambda\cdot{\bf 1})=
(s-\lambda m)B^{(m)}_{m-1}(s|\lambda\cdot{\bf 1}) =
\prod_{i=1}^m (s-\lambda i).
\label{Bernoulli_Euler_m=n_plus_1}
\ee

\section*{Acknowledgement}
The author thanks L. Fel for numerous fruitful discussions.



\end{document}